\documentclass[11pt]{amsart}
\usepackage{amsmath,amsthm,amssymb,amsfonts}
\usepackage{graphicx} 
\usepackage{epstopdf}
\makeatletter
\def\proof{\@ifstar{P\,r\,o\,o\,f}{P\,r\,o\,o\,f.\ }}
\renewcommand\th@remark{%
  \thm@headfont{\bfseries}%
  \normalfont 
  \thm@preskip\topsep \divide\thm@preskip\tw@
  \thm@postskip\thm@preskip
}
\renewenvironment{equation}{\refstepcounter{equation}$$}{\eqno{(\thesection.\theequation)}$$}
\makeatother
\newcounter{Example}[section]
\newcounter{Th}[section] \newcounter{Pr}[section] \newcounter{Lm}[section]
\newcounter{Remark}[section]
\newcounter{Def}[section]
\newcounter{lcounter}[section]
\newcounter{Corol}[section]
\newenvironment{Th}[1][\relax]
    {\medspace\refstepcounter{Th}T\,h\,e\,o\,r\,e\,m \arabic{section}.\theTh.\ \it}
    {\rm\medspace}
\newenvironment{Th.}[1][\relax]
    {\medspace\refstepcounter{Th}T\,h\,e\,o\,r\,e\,m \arabic{section}.\theTh.\ \it}
    {\rm\medspace}
\newenvironment{Pr.}[1][\relax]
    {\medspace\refstepcounter{Pr}P\,r\,o\,p\,o\,s\,i\,t\,i\,o\,n \arabic{section}.\thePr.\ \it}
    {\rm\medspace}
\newenvironment{Lm}[1][\relax]
    {\medspace\refstepcounter{Lm}L\,e\,m\,m\,a \arabic{section}.\theLm.\ \it}
    {\rm\medspace}
\newenvironment{Corol}[1][\relax]
    {\medspace\refstepcounter{Corol} C\,o\,r\,o\,l\,l\,a\,r\,y \arabic{section}.\theCorol.\ \it}
    {\rm\medspace}
\newenvironment{Remark}[1][\relax]
    {\medspace\refstepcounter{Remark}R\,e\,m\,a\,r\,k \arabic{section}.\theRemark.\rm\ }
    {\medspace}
\newenvironment{Example}[1][\relax]
    {\medspace\refstepcounter{Example}E\,x\,a\,m\,p\,l\,e \arabic{section}.\theExample.\rm\ }
    {\medspace}
\newenvironment{Def}[1][\relax]
    {\medspace\refstepcounter{Def}D\,e\,f\,i\,n\,i\,t\,i\,o\,n \arabic{section}.\theDef.\rm\ }
    {\medspace}
\newenvironment{Def.}[1][\relax]
    {\medspace\refstepcounter{Def}D\,e\,f\,i\,n\,i\,t\,i\,o\,n \arabic{section}.\theDef.\rm\ }
    {\medspace}

\def\au#1{\emph{#1}}
\def\tit#1{{#1}}
\def\R{{\mathbb R}}  
\def\H{\mathcal{H}} 
\def\cl{\mathop{\rm cl\,}}

\def\co{\mbox{\rm co}\,}
\def\cone{\mbox{\rm cone}\,}

\def\ll{\lambda}
\def\d {\partial\,}
\def\ep{\varepsilon} 
\def\Int {\mbox{\rm int\,}} 

\def\diam {\mbox{\rm diam}\,}
\def\diff {\,\frac{\,*\,}{}\,}
\def\D {\mbox{\rm strco}}
\begin{document}
\title[On Pli\'s metric on the space of strictly convex compacta]{On Pli\'s metric on the space of strictly convex compacta}
\author[M.~V.~Balashov and D.~Repov\v s]{Maxim V. Balashov and Du\v{s}an Repov\v{s}}
\address{Department of Higher Mathematics, Moscow Institute of Physics and Technology, Institutski str. 9,
Dolgoprudny, Moscow region, Russia 141700. balashov@mail.mipt.ru}
\address{Faculty of Mathematics and Physics, and Faculty of Education, University of Ljubljana, Jadranska 19, Ljubljana, Slovenia 1000.
dusan.repovs@guest.arnes.si}
\keywords{Metric space, strictly convex compactum, modulus of
convexity, set-valued mapping, strict convexity,  uniform
convexity, supporting function, Demyanov distance, Hausdorff
distance.}
\subjclass[2010]{Primary: 54A20, 52A41.  Secondary: 52A20, 52A99,
46N10.}
\begin{abstract}
We consider a certain metric on the space of all convex compacta
in $\R^{n}$, introduced by A. Pli\'s. The set of strictly convex
compacta is a complete metric subspace of the metric space of
convex compacta with respect to this metric. We present some
applications of this metric to the problems of set-valued
analysis, in particular we estimate the distance between two
compact sets with respect to this metric and to the Hausdorff
metric.
\end{abstract}
\date{\today}
\maketitle
\section{Introduction}
\def\i {\mbox{\rm int}\,}
We begin by some definitions for a finite-dimensional Euclidean
space $(\R^{n}, \|\cdot\|)$ over $\R$ with the inner product
$(\cdot,\cdot)$. Let $B_{r}(a)=\{ x\in \R^{n}\ |\ \| x-a\| \le
r\}$. Let $\cl A$ denote the \it closure \rm and $\Int A$ the \it
interior \rm of the subset $A\subset \R^{n}$. The \it diameter
\rm of the subset $A\subset \R^{n}$ is defined by $\diam A =
\sup\limits_{x,y\in A} \| x-y\|$. The \it distance \rm from the
point $x\in \R^{n}$ to the set $A\subset\R^{n}$ is given by the
formula $\varrho (x,A)=\inf\limits_{a\in A}\| x-a\|$. We shall
denote the \it convex hull \rm of a set $A\subset \R^{n}$ by $\co
A$. We shall denote the \it conic hull \rm of a set $A\subset
\R^{n}$ by $\cone A$ (cf.
\cite{Aubin,Polovinkin+Balashov,Rockafellar}).

The \it Hausdorff distance \rm between two subsets $A,B\subset
\R^{n}$ is defined as follows $h(A,B)=$
$$ = \max \left\{ \sup_{a\in A}\ \inf_{b\in B} \|a-b\| ,
\quad \sup_{b\in B}\ \inf_{a\in A} \|a-b\| \right\}=\inf\{ r>0\
|\ A\subset B+B_{r}(0),\ B\subset A+B_{r}(0) \}.
$$

The \it supporting function \rm of the subset $A\subset \R^{n}$
is defined as follows
\begin{equation}\label{**1}
s(p,A) = \sup\limits_{x\in A}(p,x),\qquad\forall p\in \R^{n}.
\end{equation}
The supporting function of any set $A$ is always lower
semicontinuous, positively uniform and convex. If the set $A$ is
bounded then the supporting function is Lipschitz continuous
\cite{Aubin,Polovinkin+Balashov}.

It follows from the Separation Theorem  (cf. \cite[Lemma
1.11.4]{Polovinkin+Balashov}) that for any convex compacta $A,B$
in $\R^{n}$
\begin{equation}\label{*3}
h(A,B) = \sup\limits_{\| p\|=1}|s(p,A)-s(p,B)|.
\end{equation}
Let $(T,\varrho)$ be a metric space. We say that a set-valued
mapping $F:(T,\varrho)\to 2^{\R^{n}}\backslash\{\emptyset \}$ is
\it upper semicontinuous \rm at the point $t=t_{0}$ if
$$
\forall \ep>0\ \exists \delta>0\ \forall t:\ \varrho
(t,t_{0})<\delta\qquad F(t)\subset F(t_{0})+B_{\ep}(0),
$$
and \it lower semicontinuous \rm at the point $t=t_{0}$ if
$$
\forall \ep>0\ \exists \delta>0\ \forall t:\ \varrho
(t,t_{0})<\delta\qquad F(t_{0})\subset F(t)+B_{\ep}(0).
$$
We say that a set-valued mapping $F:(T,\varrho)\to
2^{\R^{n}}\backslash\{\emptyset \}$ is \it continuous \rm at the
point $t=t_{0}$ if $F$ is upper and lower semicontinuous at the
point $t=t_{0}$.

We say that a set-valued mapping $F:(T,\varrho)\to
2^{\R^{n}}\backslash\{\emptyset \}$ is \it (upper, lower)
(semi)continuous \rm on the set $T$, if $F$ is (upper, lower)
(semi)continuous at any point $t_{0}\in T$.

For any convex compact set $A\subset\R^{n}$ and any vector $p\in
\R^{n}$, the subset $A(p)=\{ x\in A\ |\ (p,x)=s(p,A)\}$ is the
subdifferential of the supporting function $s(p,A)$ at the point
$p$. The set-valued mapping $\R^{n}\ni p\to A(p)$ is always upper
semicontinuous (cf. \cite{Aubin, Rockafellar}).

A convex compactum in $\R^{n}$ is called \it strictly \rm convex
if its boundary contains no nondegenerate line segments.

\begin{Def}\label{modulus} (\cite{Polyak}).
Let $E$ be a Banach space and let a subset $A\subset E$ be convex
and closed. \it The modulus of convexity \rm $\delta_{A}:\
[0,\diam A)\to [0,+\infty) $ is the function defined by
$$
\delta_{A}(\ep) = \sup\left\{ \delta\ge 0\ \left|\
B_{\delta}\left( \frac{x_{1}+x_{2}}{2}\right)\right.\subset A,\
\forall x_{1},x_{2}\in A:\ \| x_{1}-x_{2}\|=\ep \right\}.
$$\rm
\end{Def}

\begin{Def}\label{RM} (\cite{Polyak}).
Let $E$ be a Banach space and let a subset $A\subset E$ be convex
and closed. If the modulus of convexity $\delta_{A}(\ep )$ is
strictly positive for all $\ep\in (0,\diam A)$, then we call the
set $A$ \it uniformly convex \rm (\it with modulus \rm
$\delta_{A}(\cdot)$).\rm
\end{Def}

For any uniformly convex set $A$ the modulus $\delta_{A}$ is a
strictly increasing function on the segment $[0,\diam A)$.  In the
finite-dimensional case the class of strictly convex compacta
coincides with the class of uniformly convex compacta with moduli
of convexity $\delta_{A}(\ep)>0$ for all permissible $\ep>0$  (cf.
\cite{Balashov+Repovs2}).

We shall use $*$ for objects from conjugate space $E^{*}$:
$\|\cdot\|_{*}$ is the norm in $E^{*}$, $B_{1}^{*}(0)$ is the
unit closed ball in $E^{*}$ and so on.

Following Pli\'s \cite{Plis} we define the metric $\rho$
which is the main objective of the present paper.

\begin{Def} (\cite[Formula (3)]{Plis})\label{metric}
The metric $\rho$ on the space of convex compacta in $\R^{n}$ is
defined by the formula
\begin{equation}\label{metric-rho}
\rho (A,B) = \sup\limits_{\| p\|=1}h(A(p),B(p)),
\end{equation}
for any convex compacta $A,B\subset \R^{n}$.
\end{Def}

Definition 1.\ref{metric} coincides with
the definition of  the \it
Demyanov metric \rm (see its definition in
\cite[Formula (4.1)]{3metrics}) -- this was
proved in \cite{LR}.

The Hausdorff metric is the most natural metric for various
questions of set-valued analysis and its applications.
Nevertheless, there are some limitations for using this metric. For
example, if we have a sequence $\{ A_{k}\}_{k=1}^{\infty}$ of
strictly convex compact sets and $h(A_{k},A)\to 0$, then the
limit set $A$ needs not be  strictly convex. Indeed, consider
on the Euclidean plane the following ellipsoids
$$
A_{k}=\{ (x_{1},x_{2})\in\R^{2}\ |\ x_{1}^{2}+k^{2}x_{2}^{2}\le
1\}.
$$
Each set $A_{k}$ is strictly convex, but the limit set $A=\{
(x_{1},0)\in\R^{2}\ |\ x_{1}\in [-1,1]\}$ is not strictly convex.
Note, that strict convexity of the set means differentiability of
the supporting function of this set. This fact is very useful for
applications. Below we give some sufficient conditions for the
limit of a sequence of strictly convex compacta to  be
also strictly convex.

We say that a sequence of convex compacta
$\{A_{k}\}_{k=1}^{\infty}\subset \R^{n}$ is \it uniformly convex
with modulus $\delta$ \rm if $\inf\limits_{k}\diam A_{k}>0$ and
the function $\delta (\ep)$, $\ep\in [0,\inf_{k}\diam A_{k})$, is
continuous and has the property $0<\delta (\ep)\le
\delta_{A_{k}}(\ep)$ for all $\ep\in (0,\inf_{k}\diam A_{k})$ and
for all $k$.

\begin{Lm}\label{H-sc}
Let a sequence $\{A_{k}\}_{k=1}^{\infty}\subset \R^{n}$ of convex
compacta converge to a convex compactum $A$ in the Hausdorff
metric. If the sequence $\{A_{k}\}_{k=1}^{\infty}$ is uniformly
convex with modulus
$\delta$,
$\delta:(0,\ep_{0}]\to
(0,+\infty)$, then the compactum $A$ is a uniformly convex set
with the modulus $\delta_{A}(\ep)\ge \delta (\ep)$, $0<\ep\le
\ep_{0}$. In particular, this implies  strict convexity of the
set $A$.
\end{Lm}

\proof Choose arbitrary points $x,y\in A$ with $\| x-y\|<\ep_{0}$.
There are two sequences $\{x_{k}\}\subset A_{k}$,
$\{y_{k}\}\subset A_{k}$ such that $x_{k}\to x$, $y_{k}\to y$,
$k\to\infty$. For all sufficiently large $k$ we have $\|
x_{k}-y_{k}\|<\ep_{0}$. Due to the uniform convexity of the
sequence $\{ A_{k}\} $ we obtain that
$$
\frac{x_{k}+y_{k}}{2}+B_{\delta(\| x_{k}-y_{k}\|)}(0)\subset A,
$$
and
$$
\left(p, \frac{x_{k}+y_{k}}{2}\right)+\delta(\|
x_{k}-y_{k}\|)\|p\|\le s(p,A_{k}),\quad\forall p\in \R^{n}.
$$
Taking the limit $k\to\infty$, using (1.\ref{*3}) and the
continuity of the function $\delta$ we get
$$
\left(p, \frac{x+y}{2}\right)+\delta(\| x-y\|)\|p\|\le
s(p,A),\quad\forall p\in \R^{n},
$$
i.e.
$$
s\left( p,\frac{x+y}{2}+B_{\delta(\| x-y\|)}(0)\right)\le
s(p,A),\quad\forall p\in \R^{n}.
$$
By the Separation Theorem \cite{Polovinkin+Balashov,Rockafellar}
we obtain the following
$$
\frac{x+y}{2}+B_{\delta(\| x-y\|)}(0)\subset A.
$$
\qed

\section{The main properties of metric $\rho$}

In general, the subdifferential of a convex function is only upper
semicontinuous \cite{Aubin,Rockafellar}. For (not strictly) convex
compactum $A$ the sets $A(p)$ are also upper semicontinuous with
respect to $p$. This leads to the fact that in the formula
(1.\ref{metric-rho}) from Definition 1.\ref{metric} one cannot
replace sup by max.

\begin{Example}
Consider in $\R^{3}$ two sets:
$$
A=\co\left\{\{ (x_{1},x_{2},x_{3})\ |\ (x_{1}-1)^{2}+x_{2}^{2}=1;\
x_{3}=0\}\cup\{ (0,0,1)\} \right\},
$$
$$
B=\co\left\{\{ (x_{1},x_{2},x_{3})\ |\
(x_{1}-1)^{2}+x_{2}^{2}+x_{2}^{8}=1;\ x_{3}=0\}\cup\{ (0,0,1)\}
\right\}.
$$
It is easy to see that $B\subset A$, $\diam B=\diam A=\sqrt{5}$,
and $\diam A$ and $\diam B$ are attained only on the the line
segment $[(0,0,1),(2,0,0)]\subset B$.

Let $a_{k}\in \{ (x_{1},x_{2},x_{3})\ |\
(x_{1}-1)^{2}+x_{2}^{2}=1;\ x_{2}<0;\ x_{3}=0\}$ such that
$a_{k}\to (2,0,0)$. The line segment $[(0,0,1),a_{k}]$ is a
generatrix of the cone $A$ for all $k$.

Let $H_{k}$ be a supporting plane of the set $A$ such that
$[(0,0,1),a_{k}]\subset H_{k}$. Let $p_{k}$ be a unit normal
vector to the plane $H_{k}$ such that $(p_{k},a_{k})>0$. It is
easy to see that $p_{k}\to p_{0}=\frac{1}{\sqrt{5}}(1,0,2)$.

For any $k$ we have $B(p_{k})=\{ (0,0,1)\}$ and
$A(p_{k})=H_{k}\cap A = [(0,0,1),a_{k}]$. By Definition
1.\ref{metric} it follows that
$$
\rho (A,B)\ge h(A(p_{k}),B(p_{k}))=h\left( \{(0,0,1)\},\{
(0,0,1),a_{k}\}\right)=\| (0,0,1)-a_{k}\|=\sqrt{\| a_{k}\|^{2}+1},
$$
and $\sqrt{\| a_{k}\|^{2}+1}\to\sqrt{5}$, $\sqrt{\|
a_{k}\|^{2}+1}<\sqrt{5}$ for all $k$. However, the only line
segment which realizes $\diam A=\diam B=\sqrt{5}$ is the line
segment $[(0,0,1),(2,0,0)]\subset A\cap B$. Thus $\rho
(A,B)=\lim\limits_{k\to\infty}h(A(p_{k}),B(p_{k}))=\sqrt{5}$, but
for all $p$, $\| p\|=1$, $h(A(p),B(p))<\sqrt{5}$.\qed
\end{Example}

\begin{Lm}\label{lower-strict}
Let $A\subset \R^{n}$ be a convex compactum. If the set-valued
mapping $\d B_{1}(0)\ni p\to A(p)$ is lower semicontinuous, then
the compactum $A$ is strictly convex.
\end{Lm}

\proof Suppose that there exists $p\in \d B_{1}(0)$ such that the
set $A(p)$ is not a singleton. Let $\{x,y\}\subset A(p)$, $x\ne
y$, and $q=\frac{y-x}{\| y-x\|}$. Obviously, $q$ is orthogonal to
$p$.

Consider a sequence $\{ q_{k}\}_{k=1}^{\infty}\subset\cone \{
p,q\}$ such that $q_{k}\to p$, $k\to \infty$, $\| q_{k}\|=1$ and
$q_{k}\ne p$ for all $k$. Let $H_{p}^{-}=\{ z\in\R^{n}\ |\
(p,z)\le s(p, A)\}$, $H_{q_{k}}^{+}=\{ z\in\R^{n}\ |\ (q_{k},z)\ge
(y,q_{k})\}$, $H_{q}^{+}=\{ z\in\R^{n}\ |\ (q,z)\ge (y,q)\}$.

By  lower semicontinuity of $A(\cdot)$ we have for any $\ep>0$
and for all sufficiently large $k$
\begin{equation}\label{n1}
A(p)\subset A(q_{k})+B_{\ep}(0).
\end{equation}
On the other hand,
\begin{equation}\label{n2}
A(q_{k})\subset H_{q_{k}}^{+}\cap H_{p}^{-}\subset H_{q}^{+}\cap
H_{p}^{-}.
\end{equation}
Due to the inclusion (2.\ref{n2}) we obtain that
\begin{equation}\label{n3}
\varrho (x, A(q_{k}))\ge \varrho (x,H_{q}^{+}\cap H_{p}^{-})=\|
x-y\|>0.
\end{equation}
Inequality (2.\ref{n3}) implies that for all $k$
$$
x\notin A(q_{k})+\frac{\| x-y\|}{2}B_{1}(0).
$$
This contradicts the inclusion (2.\ref{n1}).\qed

\begin{Lm}\label{uniform-prezerve}
Consider a sequence $F_{k}:(T,\varrho)\to
2^{\R^{n}}\backslash\{\emptyset\}$ of set-valued mappings which
are upper (lower) semicontinuous with compact images. Let the
sequence $\{ F_{k}(t) \}_{k=1}^{\infty}$ uniformly converge to
the set-valued mapping $F:(T,\varrho)\to
2^{\R^{n}}\backslash\{\emptyset\}$ in the Hausdorff metric, i.e.
$$
\forall\ep>0\ \exists k_{\ep}\ \forall k>k_{\ep}\ \forall t\in
T\qquad h(F_{k}(t),F(t))<\ep.
$$
Then the set-valued mapping $F$ is upper (lower) semicontinuous
on $T$.
\end{Lm}

\proof The proof is a standard argument of uniform convergence.
\qed

We shall write $F_{k}\rightrightarrows F$, $t\in T$, in the case
of uniform convergence on the set $T$ of the sequence $F_{k}$ to
the mapping $F$.

\begin{Th}\label{rho-complete}  The
metric space of convex compacta in $\R^{n}$ is complete with
respect to metric $\rho$.
\end{Th}

\proof Let $\{ A_{k}\}_{k=1}^{\infty}$ be a fundamental sequence
of convex compacta with respect to metric $\rho $. This means that
$$
\forall\ep>0 \ \exists M\ \forall k,m>M\ \forall p\in\d
B_{1}(0)\qquad h(A_{m}(p),A_{k}(p))<\ep.
$$
By convexity of compact sets $A_{m}(p)$ and completeness of the
space of convex compacta with respect to the Hausdorff metric (see
\cite[Theorem 1.3.2]{Polovinkin+Balashov}) we obtain that
$A_{m}(p)\rightrightarrows A_{p}$, $p\in\d B_{1}(0)$, and the set
$A_{p}$ is convex and compact for all $p\in\R^{n}$, $\| p\|=1$.

Put
$$
A = \cl \co\bigcup\limits_{\| p\|=1}A_{p}.
$$
For any $q\in\d B_{1}(0)$ and any $x(q)\in A_{q}$ there exists a
sequence $\{ x_{m}(q)\}_{m=1}^{\infty }$ such that $x_{m}(q)\in
A_{m}(q)$ for all $m$ and $x_{m}(q)\to x(q)$. Taking a limit
$m\to\infty$ in the inequality $(p,x_{m}(p))\ge (p,x_{m}(q))$, we
obtain $(p,x(p))\ge (p,x(q))$. Hence $(p,x(p))\ge s(p,A)$ and
$x(p)\in A(p)$, i.e. $A_{p}\subset A(p)$.

The converse inclusion $A(p)\subset A_{p}$ can be proved on the
contrary with the help of separation theorem.\qed





\begin{Corol}\label{strict-complete}
The metric subspace of strictly convex compacta in $\R^{n}$ is
complete with respect to metric $\rho$.
\end{Corol}

\proof The proof is analogous to the proof of Theorem
2.\ref{rho-complete} except that all sets $A_{m}(p)$, $A_{p}$ are
singletons.\qed

Suppose that $A$, $B$ are convex compacta.  By formula
(1.\ref{*3}) we have
$$
\rho (A,B)= \sup\limits_{\| p\|=1}\sup\limits_{\| q\|=1} \left|
s(q,A(p))-s(q,B(p))\right|,
$$
 and hence
\begin{equation}\label{rho>h}
\rho (A,B)\ge\sup\limits_{\| p\|=1} \left|
s(p,A(p))-s(p,B(p))\right|= \sup\limits_{\| p\|=1} \left|
s(p,A)-s(p,B)\right|=h(A,B).
\end{equation}

Thus $\rho (A_{k},A)\to 0$ implies that $h(A_{k},A)\to 0$.

\begin{Th}\label{notLC}
The metric space of strictly convex compacta in $\R^{n}$ is not
locally compact with respect to the metric $\rho$.
\end{Th}

\proof Choose a sequence $\{A_{k}\}_{k=1}^{\infty}$ of strictly
convex compacta such that $A_{k}\subset B_{R}(0)$ for all $k$ and
there exists a nonstrictly convex compactum $A$ with
$h(A_{k},A)\to 0$.

Suppose that a subsequence $\{A_{k_{m}}\}_{m=1}^{\infty}$
converges to a compactum $B$ in the metric $\rho$.

From the estimate $\rho (A_{k_{m}},B)\ge h(A_{k_{m}},B)$ and
$h(A_{k_{m}},A)\to 0$ we get equality $B=A$. So $\rho
(A_{k_{m}},A)\to 0$. This means that
$A_{k_{m}}(p)\rightrightarrows A(p)$, $\| p\|=1$. But the set
$A_{k_{m}}(p)$ is a singleton for all $m$ and $p$. By the choice
of $A$ there exists $p_{0}\in \d B_{1}(0)$ such that the set
$A(p_{0})$ is not a singleton. Contradiction.\qed


Further we shall obtain the estimate of distance $\rho (A,B)$ via
$h(A,B)$ for some convex closed sets in a Banach space. In a
Banach space $E$ for closed convex bounded sets $A,B\subset E$ we
define
$$
\rho (A,B)=\sup\limits_{\| p\|_{*}=1}h(A(p),B(p)).
$$
If the space $E$ is reflexive then $A(p)\ne \emptyset$, $B(p)\ne
\emptyset$ for all $p\in E^{*}$.

Note that if the space $E$ contains a uniformly convex
nonsingleton set then such  space $E$ has equivalent uniformly
convex norm \cite[Theorem 2.3]{Balashov+Repovs2}. In particular,
such space $E$ is reflexive.

Note also that for any uniformly convex set $A$ we have that
$\diam A<+\infty$ and the modulus of convexity $\delta_{A}(\ep)$
is a
strictly increasing function when $\ep\in [0,\diam A)$
\cite{Balashov+Repovs2}.

\begin{Th}\label{rho<h} Let $E$ be a Banach space.
Let $A,B\subset E$ be convex closed bounded sets. Let the set $A$
be a nonsingleton and a uniformly convex set with the modulus of
convexity $\delta_{A}$. Let $\Delta=
\lim\limits_{t\to\mbox{\tiny\rm diam}\, A-0}\delta_{A}(t)$. Then
\begin{equation}\label{r<h}
\rho (A,B)\le\left\{
\begin{array}{l}
h(A,B)+\delta_{A}^{-1}(h(A,B)),\quad h(A,B)<\Delta,\\
h(A,B)\displaystyle\left(1+\frac{\diam A}{\Delta}\right), \qquad
h(A,B)\ge \Delta,
\end{array}\right.
\end{equation}
where the function $\delta_{A}^{-1}$ is the inverse function to
the function $\delta_{A}$. Furthemore, if the set $A$ is a
singleton then $\rho(A,B)=h(A,B)$.
\end{Th}

\proof Let $h=h(A,B)$.  Suppose that $A$ is not a singleton. Fix
$p\in\d B^{*}_{1}(0)$. Let $A(p)=\{a(p)\}$. Fix an arbitrary point
$b(p)\in B(p)$.

\underline{Case 1.} $h<\Delta$. Choose $t>1$ such that
$th<\Delta$.

\underline{Subcase 1.1.} $s(p,A)\ge s(p,B)$. By formula
(1.\ref{*3}) we have $0\le (p,a(p))-(p,b(p))\le h$. Let $a\in A$
be such a point that $a\in b(p)+B_{th}(0)$.

Define $ H_{A}(p)=\{ z\in E\ |\ (p,z)=s(p,A)\},\ H_{A}^{-}(p)=\{
z\in E\ |\ (p,z)\le s(p,A)\},\ H_{B}(p)=\{ z\in E\ |\
(p,z)=s(p,B)\}$.

We have $\varrho (b(p),H_{A}(p))=(p,a(p)-b(p))\le h$, $\varrho
(a,H_{A}(p)))\le \| a-b(p)\|+\varrho (b(p),H_{A}(p))\le (1+t)h$
and $A\cup B\subset H_{A}^{-}(p)$. Hence the line segment
$[a(p),a]$ belongs to the set $H_{A}(p)^{-}$. Let
$w=\frac{a(p)+a}{2}$, $\varrho (w,H_{A}(p))=\frac12\varrho
(a,H_{A}(p))\le \frac{1+t}{2}h$. By the inclusion
$$
w+\delta_{A}(\| a(p)-a\|)B_{1}(0)\subset A\subset H_{A}^{-}(p)
$$
we get
$$
\delta_{A}(\| a(p)-a\|)\le \varrho (w,H_{A}(p))\le \frac{1+t}{2}h.
$$
Hence $\| a(p)-a\|\le \delta_{A}^{-1}\left(\frac{1+t}{2}h
\right)$. Thus we obtain that
$$
\| a(p)-b(p)\|\le \| a(p)-a\|+\| a-b(p)\|\le
\delta_{A}^{-1}\left(\frac{1+t}{2}h \right)+th,
$$
i.e. $b(p)\in a(p)+\left(\delta_{A}^{-1}\left(\frac{1+t}{2}h
\right)+th\right)B_{1}(0)$. Due to the arbitrary choice of the
point $b(p)\in B(p)$ we have
$$
B(p)\subset a(p)+\left(\delta_{A}^{-1}\left(\frac{1+t}{2}h
\right)+th\right)B_{1}(0)
$$
and
$$
h(A(p),B(p))=h(\{a(p)\},B(p))\le
\delta_{A}^{-1}\left(\frac{1+t}{2}h \right)+th.
$$
Taking the limit $t\to 1+0$, we obtain that
$$
h(A(p),B(p))=h(\{a(p)\},B(p))\le \delta_{A}^{-1}\left(h \right)+h.
$$

\underline{Subcase 1.2.} $s(p,A)<s(p,B)$. Then all arguments of
the subcase 1.1 still apply except that
$$
\varrho (a, H_{A}(p))\le \varrho (a,H_{B}(p))\le \| a-b(p)\|\le
th,
$$
$\varrho (w, H_{A}(p))\le \frac{t}{2} h$, $\| a(p)-a\|\le
\delta_{A}^{-1}\left(\frac{t}{2} h \right)$. Hence
$$
h(A(p),B(p))=h(\{a(p)\},B(p))\le \delta_{A}^{-1}\left(\frac{t}{2}
h\right)+th.
$$
Taking the limit $t\to 1+0$, we obtain that
$$
h(A(p),B(p))\le \delta_{A}^{-1}\left(\frac12 h \right)+h.
$$

So again when $h<\Delta$ we have for all $p\in\d B^{*}_{1}(0)$
$$
h(A(p),B(p))\le \delta_{A}^{-1}(h)+h.
$$
Hence $\rho (A,B) = \sup\limits_{\| p\|_{*}=1}h(A(p),B(p))\le
\delta_{A}^{-1}(h)+h$.

\underline{Case 2.} $h\ge\Delta$. Then for any $t>1$ we have
$$
\rho (A,B)\le \diam A+th\le \frac{h}{\Delta}\,\diam A+th\le
h\left(t+\frac{\diam A}{\Delta}\right),\quad\forall t>1.
$$
Taking the limit $t\to 1+0$, we get
$$
\rho (A,B)\le h\left(1+\frac{\diam A}{\Delta}\right).
$$

In the case when $A$ is a singleton the equality $\rho
(A,B)=h(A,B)$ follows by definition 1.\ref{metric}.\qed

For a set $A\subset \R^{n}$, $A\subset B_{R}(a)$ for some $a\in
\R^{n}$ and $R>0$, we define \it R-strongly convex hull \rm of
the set $A$, \rm as the intersection of all closed balls of
radius $R$ each of which contains the set $A$. We shall denote the
$R$-strongly convex hull of the set $A$ by $\D_{R}A$ (cf.
\cite{Balashov+Polovinkin}).

\begin{Example}\label{best}
The estimate (2.\ref{r<h}) is exact. Consider two sets $A$ and $B$
on the Euclidean plane $\R^{2}$. Let $0<\ep<R$,
$a\left(\sqrt{2R\ep-\ep^{2}},0\right)\in \R^{2}$ and
$$
A=\D_{R}\{ B_{\ep}((0,0))\cup\{ a\}\}+B_{R}(0),\qquad B=\D_{R}\{
B_{\ep}(a)\cup\{ (0,0)\}\}+B_{R}(0).
$$
Let $p=(0,1)$. It is easy to see that $h(A,B)=\ep$, $A(p)=\left(
\ep +R\right)p$, $B(p)=\left( \sqrt{2R\ep-\ep^{2}},\ep
\right)+Rp$. Hence
$$
\rho (A,B)\ge h(A(p),B(p))=\| a\| = \sqrt{2R\ep-\ep^{2}} =
\sqrt{2R h(A,B)-h^{2}(A,B)}.
$$
The sets $A$ and $B=a-A$ are intersections of closed balls of
radius $R+\ep$ and $\delta_{A}(s)=\delta_{B}(s)\ge
(R+\ep)\delta_{\H}\left( \frac{s}{R+\ep} \right)$, where
$\delta_{\H}(s)=1-\sqrt{1-\frac{s^{2}}{4}}$ is the modulus of
convexity for the Hilbert space (see \cite[p.
63]{Lindestrauss+tzafriri}). Thus $\delta_{A}(s)=\delta_{B}(s)\ge
\frac{s^{2}}{8(R+\ep)}$, and $\delta_{A}^{-1}(t)\le
2\sqrt{2(R+\ep)t}$. So the order of $h(A,B)$ in formula
(2.\ref{r<h}) is exact.\qed
\end{Example}

\begin{Remark}
The result of Theorem 2.\ref{rho<h}  was
proved for $p$-convex sets
in  \cite[Formula (5)]{Plis}. Note that any
$p$-con\-vex set, $p>0$, in the paper \cite{Plis} is in fact the
intersection of closed balls of radius $R=\frac{1}{2p}$. From the
definition of $p$-convex set (inequality (2) of \cite{Plis}) it
follows that for any $p$-convex set $A\subset \R^{n}$, any point
$a\in \d A$ and any unit vector $w\in \{ q\in \R^{n}\ |\
(q,x-a)\le 0,\ \forall x\in A\}$ we have
$$
(w,x-a)+p\| x-a\|^{2}\le 0,\quad\forall x\in A,
$$
or
$$
A\subset B_{R}\left( a-Rw\right),\quad\mbox{\rm where}\
R=\frac{1}{2p}.
$$
Hence $A=\bigcap\limits_{\| w\|=1}B_{\frac{1}{2p}}\left(
a(w)-\frac{1}{2p}w\right)$, where $\{ a(w)\}=A(w)$.

This
also follows by results of \cite{OF}, \cite[Chapter
3]{Polovinkin+Balashov}.
\end{Remark}

\begin{Corol}\label{Cont}
Suppose that $F_{i}:(T,\varrho)\to
2^{\R^{n}}\backslash\{\emptyset\}$, $i=1,2$, are continuous (in
the metric $\rho$) set-valued mappings with strictly convex
images. Let $L:\R^{n}\to\R^{n}$ be a linear operator. Then the
set-valued mappings $F_{1}(t)+F_{2}(t)$, $LF_{1}(t)$,
$F_{2}(t)\diff F_{1}(t)=\bigcap\limits_{x\in F_{1}(t)}\left(
F_{2}(t)-x\right)$, $F_{1}(t)\cap F_{2}(t)$ (the latter two  if
nonempty) are continuous in the metric $\rho$.
\end{Corol}

\proof The proof is similar for all cases. Let us prove the
continuity of $F_{1}(t)\cap F_{2}(t)$.

The continuity of set-valued mappings $F_{i}$ in the metric $\rho$
and formula (2.\ref{rho>h}) gives the continuity of set-valued
mappings $F_{i}$ in the Hausdorff metric.

It is well known that the intersection of two continuous in the
Hausdorff metric set-valued mappings with compact strictly convex
images is also continuous in the Hausdorff metric (cf.
\cite{Aubin,Polovinkin+Balashov}). Thus the set-valued mapping
$H=F_{1}\cap F_{2}:(T,\varrho)\to
2^{\R^{n}}\backslash\{\emptyset\}$ is continuous in the Hausdorff
metric.

For any point $t=t_{0}\in T$ the set $H(t_{0})$ is a
strictly(=uniformly) convex compactum from $\R^{n}$ with some
modulus of convexity $\delta_{t_{0}}$. By Theorem 2.\ref{rho<h}
we have $\rho (H(t),H(t_{0})) \le$
$$
\le \max\left\{
h(H(t),H(t_{0}))+\delta_{t_{0}}^{-1}(h(H(t),H(t_{0})));\
\left(1+\frac{\diam
H(t_{0})}{\Delta}\right)h(H(t),H(t_{0}))\right\}\mathop{\rightarrow}\limits_{t\to
t_{0}} 0,
$$
where $\Delta = \delta_{t_{0}}(\diam H(t_{0}))$. If $H(t_{0})$ is
a singleton then $\rho (H(t),H(t_{0}))=h
(H(t),H(t_{0}))\mathop{\rightarrow}\limits_{t\to t_{0}} 0$.\qed

\section{Applications}

3.1. We prove a theorem about smooth approximation of the extremal
problem.

\begin{Th}\label{Smooth}
Let $F:(T,\varrho)\to 2^{\R^{n}}\backslash\emptyset$ be a
continuous set-valued mapping with compact convex images and
suppose that there exists $r>0$ such that for all $t\in T$
$F(t)\subset B_{r}(a(t))$ for some $a(t)\in \R^{n}$. Let $\diam
F(t)\ge d>0$ for all $t\in T$.

For any $t\in T$ and $p\in \R^{n}$, $\| p\|=1$, consider the
following problem
\begin{equation}\label{ex-p}
\max\{ (p,x)\ |\ x\in F(t)\}.
\end{equation}
Then for any $\ep\in (0,1)$ there exists an approximation
 $F_{\ep}:(T,\varrho)\to
2^{\R^{n}}\backslash\emptyset$, $F(t)\subset F_{\ep}(t)$ for all
$t\in T$, $h(F(t),F_{\ep}(t))\le \ep $ for all $t\in T$, such
that for each $t\in T$ and $p\in \R^{n}$, $\| p\|=1$, the following
problem
\begin{equation}\label{ex-p-smooth}
\max\{ (p,x)\ |\ x\in F_{\ep}(t)\}
\end{equation}
has a unique solution $F_{\ep}(t,p)=\{
f_{\ep}(t,p)\}=\arg\max\limits_{x\in F_{\ep}(t)}(p,x)$ which is
H\"older continuous with the power $\frac12$ with respect to
$h(F(t_{1}),F(t_{2}))$ for all $t_{1},t_{2}\in T$. The power
$\frac12$ is the best possible in the general case.
\end{Th}

\proof Fix $\ep\in (0,1)$. Let $R=\max\{\frac{r^{2}}{\ep},r+1\}$.
Define $F_{\ep}(t)$ as the intersection of all closed balls of
radius $R$, each of which contains the set $F(t)$. This set is
nonempty because $F(t)\subset B_{R}(a(t))$.

By \cite[formulae (5.7), (5.8)]{Balashov+Polovinkin} and
\cite[Theorem 4.4.7]{Polovinkin+Balashov} we have for all
$t_{1},t_{2}\in T$
$$
h(F_{\ep}(t_{1}),F_{\ep}(t_{2}))\le
C(\ep)h(F(t_{1}),F(t_{2})),\qquad C(\ep)=\max\left\{
\sqrt{\frac{R+r}{R-r}},\ 1+\frac{r^{2}}{R(R-r)}\right\}.
$$

By \cite[Theorem 5.4]{Balashov+Polovinkin} and \cite[Theorem
4.4.6]{Polovinkin+Balashov} we have
$$
h(F(t),F_{\ep}(t))\le \frac{r^{2}}{R}\le\ep,\qquad \forall t\in T.
$$

By the inequality $\delta_{F_{\ep}(t)}(s)\ge R\delta_{\H}\left(
\frac{s}{R} \right)$, where
$\delta_{\H}(s)=1-\sqrt{1-\frac{s^{2}}{4}}$ is the modulus of
convexity for the Hilbert space \cite[p.
63]{Lindestrauss+tzafriri}, we get
$$
\delta_{F_{\ep}(t)}(s)\ge \frac{s^{2}}{8R},\qquad \forall s\in
(0,\diam F_{\ep}(t)),
$$
and by Theorem 2.\ref{rho<h} we obtain for any $p\in\d B_{1}(0)$
  that
 $$\begin{array}{l}
\| f_{\ep}(t_{1},p) - f_{\ep}(t_{2},p)\| \le \\ \qquad\le\max\{
C(\ep)h(F(t_{1}),F(t_{2}))+\sqrt{8RC(\ep)h(F(t_{1}),F(t_{2}))};\
\left(1+\frac{2r}{\Delta}\right)h(F(t_{1}),F(t_{2}))\},
 \end{array}
 $$
where $\Delta = R\delta_{\H}\left(\frac{d}{R}\right)$. On the
other hand, we have for any convex compact set $A\subset\R^{n}$
that for some constant $C>0$ the inequality $\delta_{A}(\ep)\le
C\ep^{2}$ holds for all $\ep\in (0,\diam A)$ (see
\cite{Balashov+Repovs2}). Taking into account also Example
2.\ref{best}, we see that the power $\frac12$ is the best
possible.\qed

3.2. We consider Lipschitz selections and parametrizations of
(strictly) convex compact sets with metric $\rho$.

With any convex compact set $A\subset \R^{n}$ we can associate the
\it Steiner point \rm
$$
s(A)=\frac{1}{v_{1}}\int\limits_{\|
p\|=1}s(p,A)p\,d\mu_{n-1},\qquad v_{1}=\mu_{n}B_{1}(0),
$$
where $\mu_{n}$ is the Lebesgue measure in $\R^{n}$.

It is well known that the Steiner point is a Lipschitz selection
of convex compacta in $\R^{n}$ with the Hausdorff metric, i.e. for
any convex compacta $A,B\subset\R^{n}$ we have $s(A)\in A$ and
$$
\| s(A)-s(B)\|\le \frac{2}{\sqrt{\pi}}\frac{\Gamma\left(
\frac{n}{2}+1 \right)}{\Gamma\left( \frac{n+1}{2}\right)}\,h(A,B).
$$
The Lipschitz constant (of the order $\sqrt{n}$) above is the best
possible \cite{Przeslawski1}. See also \cite[P. 53]{Schneider},
\cite{Przeslawski}, \cite[Theorem 2.1.2]{Polovinkin+Balashov} for
details.

Using the Gauss-type formula (see \cite[formula
(2.1.15)]{Polovinkin+Balashov}, \cite[formula
(3.1)]{Przeslawski}) we obtain that
$$
\frac{1}{v_{1}}\int\limits_{\|
p\|=1}s(p,A)p\,d\mu_{n-1}=\frac{1}{v_{1}}\int\limits_{\| p\|\le
1}\nabla s(p,A)\,d\mu_{n}.
$$
Note that $\nabla s(p,A)$ exists a.e. on the ball $B_{1}(0)$.

For any convex compactum $A\subset \R^{n}$ define $U(A)=\{ p\in
B_{1}(0)\ |\ \exists \,\nabla s(p,A)\}$. The function $s(p,A)$ is
Lipschitz continuous hence $\mu_{n}U(A)=\mu_{n}B_{1}(0)$. Let
$a(A, p)=\nabla s(p,A)$ for $p\in U(A)$ and $a(A,p)=0$ for $p\in
B_{1}(0)\backslash U(A)$.

Let $A,B\subset \R^{n}$ be convex compacta and $U=U(A)\cap U(B)$,
$\mu_{n}U=\mu_{n}B_{1}(0)$. Then
$$
\| s(A)-s(B)\|\le \frac{1}{v_{1}}\int\limits_{U}\|
a(A,p)-a(B,p)\|\, d \mu_{n}\le \frac{1}{v_{1}}\int\limits_{U}\rho
(A,B)\, d \mu_{n} =\rho (A,B).
$$

Thus the Steiner point is a Lipschitz selection of convex compacta
in $\R^{n}$ with metric $\rho$ with the Lipschitz constant 1.

\def\A {{\mathbb A}}

Let $\A$ be a collection of strictly convex compacta. Then for
any $p\in\R^{n}$, $\| p\|=1$, the function $a(p)=A(p)$, $A\in\A
$, is a Lipschitz selection of the family $\A$ with the Lipschitz
constant 1 in the metric $\rho$.

\begin{Th}\label{param}
Let a collection of strictly convex compacta $\A$ from $\R^{n}$ be
uniformly bounded, i.e. there exists $M>0$ such that $\|
A\|=h(\{0\}, A)\le M$ for all $A\in\A$.

Then there exists the family of functions
\begin{equation}
f_{\ll,p}:\A\to\R^{n},\qquad (\ll,p)\in [0,1]\times\d B_{1}(0),
\end{equation}
such that for any $A\in\A$ we have
$$
A = \left\{ f_{\ll,p}(A)\ |\ \ll\in [0,1],\ p\in\d B_{1}(0)
\right\}
$$
and for any $(\ll,p)\in [0,1]\times\d B_{1}(0)$ the function
$f_{\ll,p}$ is Lipschitz on $A\in \A$ selection in the metric
$\rho$ with Lipschitz constant 1.

 Moreover, the function
$[0,1]\times\d B_{1}(0)\ni (\ll,p)\to f_{\ll,p}(A)$ is continuous
for any $A\in \A$ and the function $f_{\ll,p}$ is additive:
$f_{\ll,p}(A+B)=f_{\ll,p}(A)+f_{\ll,p}(B)$, $A,B\in\A$.
\end{Th}

\proof For any $A\in\A$ we define $f_{\ll,p}(A) = \ll a(p)+(1-\ll
)s(A)\in A$. Let $A,B\in\A$. Then (note, that $b(p)=B(p)$ for any
$p\in\d B_{1}(0)$)
$$
\begin{array}{ll}
\| f_{\ll,p}(A)-f_{\ll,p}(B)\|\le \ll\| a(p)-b(p)\|+(1-\ll )\| s(A)-s(B)\|\le\qquad\qquad\qquad\\
\qquad \qquad\qquad\qquad\qquad\qquad\le \ll\rho (A,B)+(1-\ll
)\rho (A,B)=\rho (A,B).
\end{array}
$$

Choose $\ll_{1},\ll_{2}\in [0,1]$ and $p_{1},p_{2}\in\d B_{1}(0)$
and $A\in\A$.
$$
\begin{array}{ll}
\| f_{\ll_{1},p_{1}}(A)-f_{\ll_{2},p_{2}}(A)\| = \| \ll_{1}
a(p_{1})+(1-\ll )s(A)-\ll_{2}
a(p_{2})-(1-\ll_{2} )s(A)\|\le\qquad \\
\le \| \ll_{1} a(p_{1})-\ll_{2} a(p_{2})\|+|\ll_{1}-\ll_{2}|\|
s(A)\| \le |\ll_{1}-\ll_{2}|\| a(p_{1})\|+|\ll_{2}|\| a(p_{1})-a(p_{2})\|+\\
\qquad + |\ll_{1}-\ll_{2}|\| s(A)\|\le 2|\ll_{1}-\ll_{2}|M+\|
a(p_{1})-a(p_{2})\|.
\end{array}
$$
The gradient $a(p)=\nabla s(p,A)$ for the strictly convex compact
set $A$ is uniformly continuous on the unit sphere (see
\cite[Lemma 2.2]{Balashov+Repovs2}). So the function
$[0,1]\times\d B_{1}(0)\ni (\ll,p)\to f_{\ll,p}(A)$ is uniformly
continuous.

By the Moreau-Rockafellar theorem \cite{Rockafellar} for all
$A,B\in\A$ we get $A(p)+B(p)=(A+B)(p)$ for all $p\in\d B_{1}(0)$.
Using the additive property of the Steiner point
 \cite{Polovinkin+Balashov}, \cite{Przeslawski}, \cite{Schneider}
 we obtain that $f_{\ll,p}$ is an additive
selection for all $\ll\in [0,1]$ and $\| p\|=1$.\qed

\begin{Remark}
We see from the proof of Theorem 3.\ref{param}, that the function
$$
(\ll,p)\to f_{\ll,p}(A)
$$
is uniformly continuous for any $A\in\A$. More precisely,
$f_{\ll,p}$ is Lipschitz on $\ll\in [0,1]$ (with Lipschitz
constant $2M$) and uniformly continuous on $p\in \d B_{1}(0)$.
Note that $f_{\ll,p}(A)$ is Lipschitz on $p\in \d B_{1}(0)$ if and
only if the set $A$ is an intersection of closed balls of the same
fixed radius. The last assertion follows by results of \cite{OF}
and by Theorem 4.3.2 of \cite{Polovinkin+Balashov}: \it a set $A$
is the intersection of closed balls of fixed radius $R>0$ in
Hilbert space if and only if $\| a(p)-a(q)\|\le R\| p-q\|$ for all
$p,q\in\d B_{1}(0)$. \rm Here $a(p)=A(p)$.

\end{Remark}

\section*{Acknowledgements} 
This research was supported by SRA grants P1-0292-0101,
J1-2057-0101, and BI-RU/10-11/002. The first author was supported
by  RFBR grant 10-01-00139-a, ADAP project "Development of
scientific potential of higher school" 2.1.1/11133 and projects
of FAP "Kadry" 1.2.1 grant P938 and grant 16.740.11.0128.
We thank the referee for comments and suggestions.


\bigskip

\begin{thebibliography}{99}

\bibitem{Aubin} 
\au{J.-P. Aubin, I. Ekeland,} 
\tit{Applied Nonlinear
Analysis,} 
John Wiley \& Sons Inc., New York, 1984.

\bibitem{Balashov+Polovinkin} 
\au{M. V. Balashov, E. S.
Polovinkin,} 
\tit{$M$-strongly convex subsets and their
generating sets,} 
Sbornik: Math. 191:1 (2000), 25-60.

\bibitem{Balashov+Repovs2} 
\au{M. V. Balashov, D. Repov\v{s},} 
\tit{Uniform convexity and the spliting problem for selections,}
J. Math. Anal. Appl. 360:1 (2009), 307-316.

\bibitem{3metrics} 
\au{P. Diamond, P. Kloeden, A. Rubinov, A.
Vladimirov,} 
\tit{Comparative properties of three metrics in the
space of compact convex sets,}
Set-Valued Anal. 5:3 (1997),
267-289.

\bibitem{OF}
 \au{H. Frankowska, Ch. Olech,}
 \tit{R-convexity of the integral of the set-valued functions,}
Contributions to Analysis and Geometry, John Hopkins Univ. Press,
Baltimore, Md., 1981, pp. 117-129.

\bibitem{LR} 
\au{A. Le\'sniewski and T. Rzezuchowski,}
\tit{The Demyanov Metric for Convex, Bounded Sets and Existence of
Lipschitzian Selectors,} 
J.  Convex Analysis. 18:3 (2011), in print.

\bibitem{Lindestrauss+tzafriri}
 \au{J. Lindenstrauss, L. Tzafriri,}
 \tit{Geometry of Banach Spaces - II. Functional Spaces,}
Springer-Verlag, Berlin, 1979.

\bibitem{Plis} \au{A. Plis,}
\tit{Uniqueness of optimal trajectories for non-linear control
problems,} 
Ann. Polon. Math. 29 (1975), 397-401.

\bibitem{Polovinkin+Balashov}
\au{E. S. Polovinkin, M. V.
Balashov,}
\tit{Elements of Convex and Strongly Convex Analysis,}
Fizmatlit, Moscow, 2007. (in Russian).

\bibitem{Polyak}
\au{B. T. Polyak,} 
\tit{Existence theorems and convergence of minimizing sequences
in extremum problems  with restrictions,} 
Soviet Math. 7 (1966),
72-75.

\bibitem{Przeslawski1} 
\au{K. Przeslawski,} 
\tit{Linear and Lipschitz continuous selectors for the family of
convex sets in Euclidean vector spaces,}
Bull. Polish Acad. Sci.
Math. 33:1-2 (1985), 31--34.

\bibitem{Przeslawski} 
\au{K. Przeslawski,} 
\tit{Lipschitz continuous selectors, Part I: Linear
selectors,} 
J. Convex Anal. 5:2 (1998), 249-267.

\bibitem{Rockafellar} 
\au{R.T. Rockafellar,} 
\tit{Convex Analysis,} 
Princeton University Press,
Princeton, NJ, 1970.

\bibitem{Schneider} 
\au{R. Schneider,} 
\tit{Convex Bodies: The Brunn-Minkowski
Theory,} 
Cambridge Univ. Press, 1993.
\end{thebibliography}
\end{document}